\documentclass[french,10pt]{article} 
\usepackage[latin1]{inputenc}
\usepackage{fontenc}
\usepackage{amsmath}
\usepackage{amssymb}
\usepackage{graphicx}
\usepackage{epsfig}
\usepackage{babel}
\usepackage[all]{xy}
\newtheorem{lemme}{Lemme}

\newtheorem{proposition}{Proposition}

\def\R{{\mathbb R}} 
\def\Z{{\mathbb Z}} 

\def\N{{\mathbb N}}

\def\V{\text{Vol}}
\def\s{\text{sys}}

\def\l{\text{long}}

\begin{document}

\title{Invariant d'Hermite des jacobiennes de graphes pondérés}

\author{Florent BALACHEFF\footnote{UMR 5149, Institut de Math\'ematiques et de Modélisation
de  Montpellier, Universit\'e Montpellier II  Case Courrier 051 - Place Eug\`ene Bataillon 34095
 Montpellier CEDEX 5, France e-mail: balachef@math.univ-montp2.fr}} 

\date{}

\maketitle

\bigskip

\begin{abstract}

A tout graphe pondéré de premier nombre de Betti $b$ est naturellement associé un réseau de dimension $b$, défini de manière analogue que la jacobienne pour une surface de Riemann. Cette classe de réseaux générée par des graphes est particulièrement intéressante. Nous démontrons ici une majoration de l'invariant d'Hermite d'un tel réseau en fonction de $b$ dont l'ordre de grandeur est $\ln b$. Cet ordre de grandeur est optimal : il est réalisé par l'invariant d'Hermite de la jacobienne d'un graphe systoliquement économique.
\end{abstract}

\bigskip


\noindent {\it Mots clefs :} Constante systolique, complexité, graphe, invariant d'Hermite, jacobienne, systole.
\bigskip

\section*{Introduction}

Soit $n$ un entier non nul. L'étude de la densité des réseaux d'un espace euclidien $(\R^n,<.,.>)$ est un sujet classique et prend la forme suivante :  étant donné un réseau $\Lambda \subset \R^n$, nous définissons le déterminant de $\Lambda$, noté $\det(\Lambda)$, comme le carré du volume euclidien du domaine fondamental du réseau et sa norme minimale par la formule
$$
||\Lambda||= \min \{<\lambda,\lambda> \mid \lambda \in \Lambda \setminus \{0\}\}.
$$
 L'invariant d'Hermite du réseau est la quantité
$$
\mu(\Lambda)=\frac{||\Lambda||}{\sqrt[n]{\det(\Lambda)}}
$$
et code la densité du réseau. La densité maximale en dimension $n$ correspond à la {\it constante d´Hermite} :
\begin{equation}  \label{herm}
\gamma_n=\sup\{ \mu(\Lambda) \mid \Lambda \text{ réseau de } \R^n\}. 
\end{equation}
Cette quantité est bien définie et vérifie l'encadrement suivant (voir \cite{conslo}) :
\begin{equation} \label{herm2}
\frac{n}{2\pi e} \lesssim \gamma_n \lesssim \frac{1.744 n}{2\pi e}.
\end{equation}

Dans la définition (\ref{herm}), on peut considérer la borne supérieure non plus sur tous les réseaux mais sur un sous-ensemble de réseaux de $\R^n$. P. Buser et P. Sarnak ont étudié dans \cite{busesarn94} la borne supérieure des invariants d'Hermite des réseaux symplectiques et ont montré qu'elle vérifiait l'inégalité inférieure dans la formule (\ref{herm2}). Ils ont également montré le résultat suivant : la borne supérieure des invariants d'Hermite sur l'ensemble des réseaux formé des jacobiennes de surfaces de Riemann de genre $g$, que l'on notera $\eta_{2g}$, vérifie 
\begin{equation} \label{herm3}
c \ln g \lesssim \eta_{2g} \lesssim \frac{3}{\pi} \ln (4g+3),
\end{equation}
où $c$ est une constante positive.

\medskip

Le but de cet article est de prouver un résultat analogue pour la jacobienne des graphes pondérés. Rappelons tout d'abord quelques définitions. Un {\it graphe} $\Gamma=(V,E)$ est un complexe simplicial de dimension 1. C'est la donn\'ee  d'une paire d'ensembles $(V,E)$, o\`u $V$ d\'esigne les sommets et $E$ les ar\^etes. La {\it valence} d'un sommet est le nombre d'arêtes incidentes à ce sommet et un graphe sera dit {\it $k$-régulier}, pour $k\in \N^\ast$, si la valence de chacun de ses sommets est constante égale à $k$. 
Dans ce qui suit, les graphes seront supposés  connexes, finis et leurs sommets de valence au moins $2$.  Un {\it graphe pond\'er\'e} est une paire $(\Gamma,w)$ o\`u $\Gamma=(V,E)$ est un graphe et $w$ est une {\it fonction poids} sur les ar\^etes $w:E \rightarrow \mathbb R_+$. Nous dirons d'un graphe pondéré $(\Gamma,w)$ qu'il est {\it combinatoire} si sa fonction poids est constante égale à $1$, et nous le noterons simplement $\Gamma$. Le type d'homotopie d'un graphe donné $\Gamma=(V,E)$ est caractérisé par le nombre $b_1(\Gamma)$ de cycles indépendants, appelé {\it premier nombre de Betti} ou {\it nombre cyclomatique}. On a la formule $b_1(\Gamma)=|E|-|V|+1$ où $|X|$ désigne le cardinal d'un ensemble fini $X$. Notons qu'à premier nombre de Betti fixé, les graphes considérés sont en nombre fini à homéomorphisme près.

Etant donné un graphe pondéré $(\Gamma,w)$ de premier nombre de Betti $b\geq 1$, nous choisissons pour chaque arête une orientation arbitraire et nous noterons ${\mathbb E}=\{e_i\}_{i=1}^k$ l'ensemble de ces arêtes orientées. Soit 
$$
{\cal{C}}(\Gamma,\R) =\{ \sum_{i=1}^k a_i.e_i \mid  a_i \in \R \text{ pour } i=1, \ldots, k\}
$$
l'espace vectoriel engendré par les arêtes orientées. Cet espace co\"incide avec l'espace des chaînes simpliciales du complexe simplicial $\Gamma$. Il est muni du produit scalaire naturel $<e_i,e_j>_w=w(e_i)\delta_{ij}$ pour $1\leq i \leq j\leq k$, où $\delta_{ij}$ désigne le symbole de Kronecker (voir \cite{bacharnag97}, p. 191). L'homologie de $\Gamma$ de dimension $1$ à coefficients réels $H_1(\Gamma,\R)$ est plongée dans ${\cal{C}}(\Gamma,\R)$ comme un sous-espace vectoriel de dimension $b$ et on note encore $<.,.>_w$ la restriction du produit scalaire à ce sous-espace. L'homologie de $\Gamma$ de dimension $1$ à coefficients entiers $H_1(\Gamma,\Z)$, en l'absence de torsion dans ce cadre unidimensionnel, constitue un réseau du sous-espace 
$H_1(\Gamma,\R)$ (comparer avec \cite{bacharnag97}). 

On note $\Lambda(\Gamma,w)= H_1(\Gamma,\Z)\subset (H_1(\Gamma,\R),<.,.>_w)$ le réseau de dimension $b$ ainsi défini et on l'appelle {\it jacobienne} du graphe pondéré $(\Gamma,w)$ : il est clair que l'invariant d'Hermite de ce réseau ne dépend pas des orientations d'arêtes initialement choisies. Posons
$$
\rho_b=\sup \{ \mu(\Lambda(\Gamma,w)) \mid (\Gamma,w) \text{ graphe pondéré de premier nombre de Betti } b \}.
$$

Notre résultat principal s'énonce alors de la manière suivante.
\medskip

\noindent {\bf Théorème} {\it Pour $b$ suffisamment grand,  
\begin{equation} \label{theo1}
\frac{1}{6e}\log_2 b \lesssim \rho_b \lesssim 4\log_2 (\frac{8}{3}b),
\end{equation}
où $\log_2$ désigne le logarithme en base $2$.
}

\medskip

Nous pouvons, sous certaines restrictions, améliorer ces deux inégalités.

\noindent 1) Pour une infinité de valeurs $\{b_m\}_m$, il existe un graphe combinatoire \linebreak $3$-réguliers $G_m$ de premier nombre de Betti $b_m$ (construit dans \cite{marg88}) pour lequel (voir inégalité (\ref{amelio2})) 
$$
\mu(\Lambda(G_m))\gtrsim \frac{4}{9e} \log_2 b_m,
$$

\noindent 2) Tout graphe combinatoire $\Gamma$, dont la valence en chaque sommet est au moins $3$, vérifie (voir inégalité (\ref{amelio1}))
$$
\mu(\Lambda(\Gamma)) \lesssim 2\log_2 b.
$$

\bigskip

La suite de cet article est consacré aux démonstrations de ces résultats. Nous allons montrer que l'étude de l'invariant d'Hermite de la jacobienne d'un graphe pondéré est équivalente à l'étude du problème combinatoire suivant : étant donné un graphe, borner la complexité de ce graphe par sa systole (voir la section $1$ pour les définitions de ces quantités). Dans la première partie de ce papier, nous donnons la démonstration de l'inégalité supérieure dans la formule (\ref{theo1}). Dans la seconde partie, nous majorons pour tout graphe sa complexité par son volume unidimensionnel, et obtenons ainsi  à l'aide de graphes systoliquement économiques - graphes dont le rapport volume sur systole est suffisamment petit - l'estimée inférieure annoncée pour $\rho_b$. Nous démontrons pour finir les améliorations 1) et 2).

\bigskip

\section{Complexité et systole d'un graphe}

Nous allons montrer la proposition suivante, qui démontrera la majoration de $\rho_b$ annoncée dans le théorème par un calcul d'équivalent élémentaire.

\begin{proposition}
Pour tout graphe pondéré $(\Gamma,w)$ de premier nombre de Betti $b\geq 2$,
\begin{equation} \label{sup}
\mu(\Lambda(\Gamma,w))\leq 4 \left (\prod_{k=2}^b\log_2 (\frac{8}{3} k)\right )^{1/b}.
\end{equation}
\end{proposition}

\medskip

\noindent {\bf Démonstration.} Commen\c cons par réduire le problème. On peut construire facilement un graphe $3$-régulier $\Gamma'=(V',E')$  et une fonction poids $w'$ sur $\Gamma'$ tels que $(\Gamma,w)$ soit obtenu à partir de $(\Gamma',w')$ en contractant les arêtes de poids nul en un point. Notons que $\mu(\Lambda(\Gamma,w))=\mu(\Lambda(\Gamma',w'))$. Comme l'application 
\begin{eqnarray} 
\nonumber  \mu_{\Gamma'} : \R^{|E'|} & \rightarrow &  \R_+ \\
\nonumber w & \rightarrow & \mu(\Lambda(\Gamma',w))
\end{eqnarray}
est continue et invariante par composition avec les dilatations, on en déduit que pour tout $\epsilon>0$, il existe une fonction poids $w_\epsilon$ sur $\Gamma'$ telle que :

\noindent  - Pour toute arête $e \in E', w_\epsilon(e) \in \N$,

\noindent  - $| \mu(\Lambda(\Gamma,w))-\mu(\Lambda(\Gamma',w_\epsilon))|< \epsilon$.

Il suffit donc de démontrer le résultat annoncé pour un graphe combinatoire ($w=1$) dont les sommets sont de valence $2$ ou $3$. 

\medskip

Etant donné un graphe combinatoire $\Gamma$, nous pouvons exprimer l'invariant d'Hermite de sa jacobienne en fonction des quantités suivantes :

\noindent - La {\it systole}, qui dans le cas d'un graphe pondéré $(\Gamma,w)$ est définie comme la plus petite longueur d'un circuit simple de $\Gamma$ et est notée $\s(\Gamma,w)$. Dans le cas d'un graphe $\Gamma$ combinatoire, nous noterons simplement $\s(\Gamma)$ cette quantité. 

\noindent - La {\it complexité}, qui est définie comme le nombre d'arbres maximaux du graphe $\Gamma$ (voir \cite{biggs}) et est notée $\kappa(\Gamma)$. 

D'après \cite{bacharnag97}, étant donné un graphe combinatoire $\Gamma$,
\begin{center}
$||\Lambda(\Gamma)||=\s(\Gamma)$ et $\det(\Lambda(\Gamma))=\kappa(\Gamma)$.
\end{center}
\noindent On en déduit l'égalité
$$
\mu(\Lambda(\Gamma))=\frac{\s(\Gamma)}{\sqrt[b]{\kappa(\Gamma)}}.
$$

\bigskip

Supposons pour la suite de cette démonstration que $\Gamma$ est un graphe combinatoire dont les sommets sont de valence $2$ ou $3$. Notre majoration revient donc à estimer inférieurement la complexité d'un tel graphe par sa systole. Nous obtenons le résultat suivant, qui implique l'estimée (\ref{sup}) :
\begin{lemme}
\begin{equation} \label{estsup}
\kappa(\Gamma)\geq \frac{\s(\Gamma)^b}{4^b \prod_{k=2}^b \log_2 (\frac{8}{3} k)}. 
\end{equation}
\end{lemme}

\medskip

\noindent {\bf Démonstration du lemme.} Tout d'abord, supposons que $b=2$. Comme les sommets de $\Gamma$ sont de valence $2$ ou $3$, on a deux classes d'homéomorphismes possibles pour $\Gamma$ : la classe $8_1$ et la classe $8_2$ (voir figure 1).

\begin{figure}[h]
\begin{center}
\includegraphics[width=10cm]{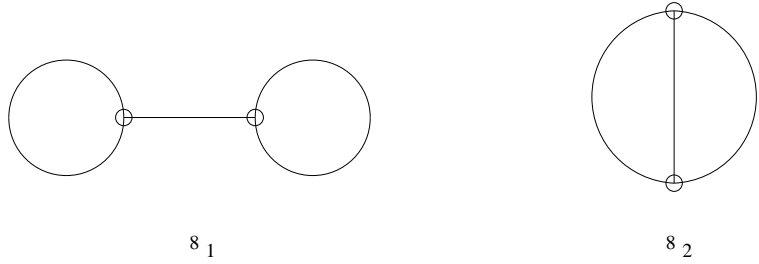}
\caption{Les graphes $8_1$ et $8_2$}
\end{center}
\end{figure}

Il est clair que $\kappa(8_1) \geq (\s(8_1))^2$ et $\kappa(8_2) \geq (\s(8_2)/2)^2$, d'où (\ref{estsup}) dans ce cas.

\medskip

On suppose maintenant $b>2$. Si 
$$
\left[\frac{\s(\Gamma)}{2 \log_2 (\frac{8}{3} b)}\right]=0,
$$
où $[n]$ désigne la partie entière d'un entier $n$, le résultat est immédiat comme la complexité d'un graphe est un entier non nul. Supposons donc 
$$
\left[\frac{\s(\Gamma)}{2 \log_2 (\frac{8}{3} b) }\right]\geq 1.
$$

On note $\tilde{\Gamma}$ le graphe défini à partir de $\Gamma$ de la manière suivante. Etant donné un sommet $v$ de valence $2$, on considère le graphe obtenu en supprimant le sommet $v$ et les deux arêtes $e_1$ et $e_2$ incidentes à ce sommet, et en ajoutant une nouvelle arête reliant les deux sommets de $e_1$ et $e_2$ restants (voir figure $2$). 
\begin{figure}[h]
\begin{center}
\includegraphics[width=4cm]{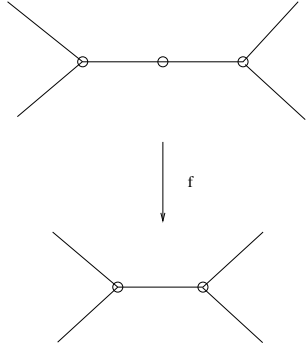}
\caption{Elimination d'un sommet de valence $2$}
\end{center}
\end{figure}
On répète l'opération pour tout sommet de valence $2$ et on obtient ainsi le graphe $\tilde{\Gamma}$. C'est un graphe $3$-régulier et on note $f : \Gamma \rightarrow \tilde{\Gamma}$ l'application topologique naturelle qui envoie uniformément une suite maximale d'arêtes adjacentes de $\Gamma$ dont les sommets intermédiaires sont de valence $2$ sur l'arête correspondante de $\tilde{\Gamma}$. Soit $\tilde{\gamma}$ une courbe réalisant sa systole et $\gamma$ son image inverse par $f$. La systole est estimée supérieurement pour les graphes $3$-réguliers de premier nombre de Betti $b$ de la manière suivante (voir \cite{erdosach63}) : 
$$
\s(\tilde{\Gamma}) \leq [2 \log_2 (\frac{2}{3} b) +3].
$$
 Donc il existe une arête $e \in \tilde{\gamma} \subset \tilde{\Gamma}$ telle que, si $C=f^{-1}(e)$, on ait 
$$
\l(C) \geq \left[\frac{\s(\Gamma)}{2 \log_2 (\frac{8}{3} b) }\right].
$$
En effet, sinon, $\l(\gamma) < \s(\Gamma)$, d'où une contradiction. 

Soit $\Gamma \setminus C$ le complémentaire de la suite d'arêtes adjacentes $C$ dans $\Gamma$ : c'est un sous-graphe de $\Gamma$ vérifiant $b_1(\Gamma \setminus C)=b_1(\Gamma)-1$. Tout arbre maximal $T$ du graphe $\Gamma \setminus C$ fournit de manière évidente au moins $[\s(\Gamma)/(2 \log_2 (\frac{8}{3} b) )]$ arbres maximaux $T'$ de $\Gamma$, et tous les arbres maximaux ainsi obtenus, lorsque l'on fait varier l'arbre initial $T$ parmi les arbres maximaux de $\Gamma \setminus C$,  sont deux à deux distincts. On en déduit :
$$
\kappa(\Gamma) \geq \kappa(\Gamma\setminus C).\left [\frac{\s(\Gamma)}{2 \log_2 (\frac{8}{3} b) }\right]\geq \kappa(\Gamma\setminus C).\frac{\s(\Gamma)}{2( 2\log_2 (\frac{8}{3} b))}.
$$

Comme $\s(\Gamma \setminus C) \geq \s(\Gamma)$, on obtient le résultat par récurrence. \hspace{\stretch{1}} $\Box$\\

\noindent {\it Remarque.} Nous avons donc relié la complexité et la systole de tout graphe $\Gamma$ de premier nombre de Betti $b\geq 2$ de la manière suivante :
$$
\sqrt[b]{\kappa(\Gamma)}\geq \frac{\s(\Gamma)}{4\left (\prod_{k=2}^b(\log_2 (\frac{8}{3} k) )\right )^{1/b}}\geq \frac{\s(\Gamma)}{4\log_2 (\frac{4}{3} b)}. 
$$
De l'inégalité (\ref{estinf}), nous déduisons que l'ordre de grandeur dans cette estimée inférieure est optimale.

\medskip


\section{Complexité et volume d'un graphe} 

Nous allons dans cette section définir la notion de constante systolique pour un graphe pondéré, et ainsi pouvoir introduire une famille de graphes systoliquement éco\-nomiques. Cette famille nous permet, à l'aide d'une majoration de la complexité d'un graphe par son volume, de démontrer l'inégalité inférieure annoncée dans la formule (\ref{theo1}).

\medskip

Etant donné un graphe pondéré $(\Gamma,w)$ de premier nombre de Betti $b\geq 2$, le problème systolique peut être formulé comme suit. On définit la {\it constante systolique} de $\Gamma$ par
$$
\sigma(\Gamma)=\inf_w \frac{\V(\Gamma,w)}{\s(\Gamma,w)}
$$
où $\V(\Gamma,w)$ désigne le $1$-volume de $\Gamma$ (la somme des poids de ses arêtes) et où l'infimum est pris sur l'ensemble des fonctions poids de $\Gamma$. Nous avons, d'après \cite{bolsze02}, 
\begin{equation} \label{syst}
\sigma(\Gamma)\geq \frac{3}{2}\frac{b-1}{\log_2 (b-1)+\log_2 \log_2 (b-1)+4}. 
\end{equation}
D'autre part, comme il a été démontré dans  \cite{babebala04}, il existe pour chaque $b\geq 2$ un graphe combinatoire de premier nombre de Betti $b$ que nous noterons $\Gamma_b$ et que nous appelerons {\it systoliquement économique} vérifiant asymptotiquement
$$ 
\frac{\V(\Gamma_b)}{\s(\Gamma_b)} \lesssim 6 \frac{b}{\log_2 b},
$$ 
où $\V(\Gamma_b)$ est le $1$-volume de $\Gamma_b$ pour la fonction poids constante égale à $1$.

\medskip

Nous pouvons minorer l'invariant d'Hermite de la jacobienne d'un graphe pondéré par sa constante systolique comme suit :
\begin{proposition}
Pour tout graphe pondéré $(\Gamma,w)$ de premier nombre de Betti $b \geq 1$, on a 
$$
\mu(\Lambda(\Gamma,w))\geq \frac{\sqrt[b]{b!}}{\sigma(\Gamma,w)}.
$$
\end{proposition}

\medskip

\noindent {\bf Démonstration.} De même que dans le début de la section 1, nous nous ramenons à un graphe combinatoire. On estime supérieurement la complexité de ce graphe par son volume de la manière suivante. Tout arbre maximal $T$ de $\Gamma$ est entièrement déterminé par les $b$ arêtes de  $\Gamma \setminus T$, d'où
$$
\kappa(\Gamma) \leq C_{\V(\Gamma)}^b \leq \frac{\V(\Gamma)^b}{b!}.
$$
On en déduit immédiatement le résultat. \hspace{\stretch{1}} $\Box$\\

Les graphes $\Gamma_b$ vérifient donc
$$
\mu(\Lambda(\Gamma_b))\gtrsim \frac{1}{6}\frac{ \sqrt[b]{b!}}{b} \log_2 b,
$$
d'où la minoration annoncée en (\ref{theo1}) par application de la formule de Stirling :
\begin{equation} \label{estinf}
\mu(\Lambda(\Gamma_b))\gtrsim \frac{1}{6e} \log_2 b. 
\end{equation}

\bigskip

Nous expliquons maintenant les améliorations 1) et 2) annoncées dans l'introduction.

\medskip

\noindent 1) Nous pouvons, en utilisant les familles de graphes à grand tour de taille mises en évidence dans \cite{marg88}, améliorer pour certaines valeurs de $b$ l'estimée (\ref{estinf}).  \linebreak G.A. Margulis a construit, pour une famille infinie $\{b_m\}_m$ de valeurs, une famille $\{G_m\}_m$ de graphes $3$-réguliers de premier nombre de Betti $\{b_m\}_m$ pour lesquels 
$$
\s(G_m) \gtrsim \frac{4}{3} \log_2 b_m.
$$
Nous en déduisons
\begin{equation} \label{amelio2}
\mu(\Lambda(G_m))\gtrsim \frac{4}{9e} \log_2 b_m.
\end{equation}

\medskip

2) Si l'on se restreint à certains graphes, nous pouvons estimer supérieurement l'invariant d'Hermite de la jacobienne par la constante systolique :
\begin{proposition}
Pour tout graphe combinatoire $\Gamma$ dont chaque sommet est de valence au moins trois et de premier nombre de Betti $b\geq 2$, 
$$
\mu(\Lambda(\Gamma))\leq \frac{3(b-1)}{\sigma(\Gamma)}.
$$
\end{proposition}

\noindent {\bf Démonstration.} Comme la valence en chaque sommet du graphe est au moins $3$, on a l'inégalité 
$$
\V(\Gamma) \leq 3(b-1),
$$
donc
$$
\mu(\Lambda(\Gamma)) =\frac{\s(\Gamma)}{\sqrt[b]{\kappa(\Gamma)}} \leq  \s(\Gamma) \leq  \V(\Gamma)\frac{\s(\Gamma)}{\V(\Gamma)} \leq  \frac{3(b-1)}{\sigma(\Gamma)}.
$$
 \hspace{\stretch{1}} $\Box$\\

On en déduit une amélioration de l'inégalité (\ref{sup}) pour cette classe de graphes, en vertu de l'estimée systolique (\ref{syst}) : pour tout graphe $\Gamma$ dont chaque sommet est de valence au moins trois et de premier nombre de Betti $b\geq 2$, 
\begin{equation} \label{amelio1}
\mu(\Lambda(\Gamma))\leq 2(\log_2(b-1)+\log_2\log_2(b-1)+4) \lesssim 2\log_2 b.
\end{equation}

\bigskip

\noindent {\bf Remerciements.} Je remercie Roland Bacher de m'avoir suggérer cette étude, ainsi que Ivan Babenko pour les nombreuses discussions que nous avons partagé.

\end{document}